\documentclass[12pt]{article}
\usepackage{subeqn}
\usepackage{graphicx}
\usepackage{ifpdf}
\ifpdf  \fi
\usepackage{pstricks}
\usepackage{amsfonts,amssymb,slashbox}
\usepackage{mathptmx,helvet,courier,makeidx,multicol,footmisc}
\usepackage[numbers]{natbib}
\bibpunct{(}{)}{;}{a}{,}{,}
\usepackage[bookmarksnumbered=true, pdfauthor={Wen-Long Jin}]{hyperref}

\newcommand{\commentout}[1]{}

\newcommand{\ba}{\begin{array}}
        \newcommand{\ea}{\end{array}}
\newcommand{\bc}{\begin{center}}
        \newcommand{\ec}{\end{center}}
\newcommand{\bdm}{\begin{displaymath}}
        \newcommand{\edm}{\end{displaymath}}
\newcommand{\bds} {\begin{description}}
        \newcommand{\eds} {\end{description}}
\newcommand{\ben}{\begin{enumerate}}
        \newcommand{\een}{\end{enumerate}}
\newcommand{\beq}{\begin{equation}}
        \newcommand{\eeq}{\end{equation}}
\newcommand{\bfg} {\begin{figure}[h]}
        \newcommand{\efg} {\end{figure}}
\newcommand{\bi} {\begin {itemize}}
        \newcommand{\ei} {\end {itemize}}
\newcommand{\bqn}{\begin{eqnarray}}
        \newcommand{\eqn}{\end{eqnarray}}
\newcommand{\bqs}{\begin{eqnarray*}}
        \newcommand{\eqs}{\end{eqnarray*}}
\newcommand{\bsl} {\begin{slide}[8.8in,6.7in]}
        \newcommand{\esl} {\end{slide}}
\newcommand{\bss} {\begin{slide*}[9.3in,6.7in]}
        \newcommand{\ess} {\end{slide*}}
\newcommand{\btb} {\begin {table}}
        \newcommand{\etb} {\end {table}}
\newcommand{\bvb}{\begin{verbatim}}
        \newcommand{\evb}{\end{verbatim}}

\newcommand{\m}{\mbox}

\newcommand {\pd}[2] {{\frac {\partial {#1}} {\partial {#2}}}}

\newcommand{\cas}[1]{{{\left \{ \ba #1 \ea \right. }}}

\newcommand{\reff}[1] {{{Figure \ref {#1}}}}
\newcommand{\refe}[1] {{(\ref {#1})}}

\def\a          {{\alpha}}

\def\pmb#1{\setbox0=\hbox{$#1$}%
   \kern-.025em\copy0\kern-\wd0
   \kern.05em\copy0\kern-\wd0
   \kern-.025em\raise.0433em\box0 }

\def\eop{{\hfill $\blacksquare$}}
\def\r{{\rho}}
\newtheorem{theorem}{Theorem}[section]

\newtheorem{corollary}[theorem]{Corollary}

\def\dx     {{\Delta x}}
\def\dt     {{\Delta t}}

\begin{document}
\title{Continuous kinematic wave models of merging traffic flow}
\author{Wen-Long Jin \footnote{Department of Civil and Environmental Engineering, California Institute for Telecommunications and Information Technology, Institute of Transportation Studies, 4000 Anteater Instruction and Research Bldg, University of California, Irvine, CA 92697-3600. Tel: 949-824-1672. Fax: 949-824-8385. Email: wjin@uci.edu. Corresponding author}}
\maketitle
\begin{abstract}
Merging junctions are important network bottlenecks, and a better understanding of merging traffic dynamics has both theoretical and practical implications. In this paper, we present continuous kinematic wave models of merging traffic flow which are consistent with discrete Cell Transmission Models with various distribution schemes. In particular, we develop a systematic approach to constructing kinematic wave solutions to the Riemann problem of merging traffic flow in supply-demand space. In the new framework, Riemann solutions on a link consist of an interior state and a stationary state, subject to admissible conditions such that there are no positive and negative kinematic waves on the upstream and downstream links respectively. In addition, various distribution schemes in Cell Transmission Models are considered entropy conditions. In the proposed analytical framework, we prove that the stationary states and boundary fluxes exist and are unique for the Riemann problem for both fair and constant distribution schemes. We also discuss two types of invariant merge models, in which local and discrete boundary fluxes are the same as global and continuous ones. With numerical examples, we demonstrate the validity of the analytical solutions of interior states, stationary states, and corresponding kinematic waves. Discussions and future studies are presented in the conclusion section.
\end{abstract}

\section{{Introduction}}

In a road network, bottlenecks around merging, diverging, and other network junctions can cause the formation, propagation, and dissipation of traffic queues. For example, upstream queues can form at a merging junction due to the limited capacity of the downstream branch; at a diverging junction, when a queue forms on one downstream branch, the flow to the other downstream branch will be reduced, since the upstream link will be blocked  due to the First-In-First-Out principle on the upstream link \citep{papageorgiou1990assignment}. Theoretically, a better understanding of merging traffic flow will be helpful for analyzing and understanding traffic dynamics on a road network \citep{nie2008oscillatory,jin2009network}.
Practically, it will be helpful for designing metering schemes \citep{papageorgiou2002freeway} and for understanding dynamic user equilibrium and drivers' route choice behaviors \citep{Peeta2001dta}.

In the literature, a number of microscopic models have been proposed for vehicles' merging behaviors \citep[e.g.][]{hidas2002changing,hidas2005modelling}. Due to the complicated interactions among vehicles, microscopic models are not suitable for analyzing merging traffic dynamics at the network level. At the macroscopic level, many attempts have been made to model merging traffic flow in the line of the LWR model \citep{lighthill1955lwr,richards1956lwr}, which describes traffic dynamics with kinematic waves, including shock and rarefaction waves, in density ($\r$), speed ($v$), and flux ($q$). Based on a continuous version of traffic conservation, $\pd{\r}{t}+\pd{q}x=0$, and an assumption of a speed-density relationship, $v=V(\r)$, the LWR model can be written as
\bqn
\pd{}{t}\r+\pd {}{x} \r V(\r)=0, \label{lwr}
\eqn
which is for a homogeneous road link with time and location independent traffic characteristics, such as free flow speed, jam density, capacity, and so on. In general, $V(\r)$ is a non-increasing function, and $v_f=V(0)$ is the free flow speed. In addition, $q=Q(\r)\equiv \r V(\r)$ is unimodal with capacity $C=Q(\r_c)$, where $\r_c$ is the critical density. If we denote the jam density by $\r_j$, then $\r\in[0,\r_j]$.

In one line, \citet{daganzo1995ctm} and \citet{lebacque1996godunov} extended the Godunov discrete form of the LWR model and developed a new framework, referred to as Cell Transmission Models (CTM) hereafter, for numerically computing traffic flows through merging, diverging, and general junctions. In this framework, local traffic demand (i.e., sending flow) and supply (i.e., receiving flow) functions are introduced, and boundary fluxes through various types of junctions can be written as functions of upstream demands and downstream supplies. In the CTM framework, various distribution schemes can be employed to uniquely determine merging flows from all upstream links \citep{daganzo1995ctm,lebacque1996godunov,jin2003merge,ni2005merge}. These discrete merge models are physically intuitive, and stationary distribution patterns of merging flows have been observed and calibrated for various merging junctions \citep{cassidy2005driver,bargera2009_merge}. Such merge models are discrete in nature and can be used to simulate macroscopic traffic dynamics efficiently. However, there have been no systematic approaches to analyzing kinematic waves arising at a merging junction with these models.

In the other line, \citet{holden1995unidirection} and \citet{coclite2005network} attempted to solve a Riemann problem of a highway intersection with $m$ upstream links and $n$ downstream links. 
In both of these studies, all links are homogeneous and have the same speed-density relations, and traffic dynamics on each link are described by the LWR model. In \citep{holden1995unidirection}, the Riemann problem is solved by introducing an entropy condition that maximizes an objective function in all boundary fluxes. In \citep{coclite2005network}, the Riemann problem is solved to maximize total flux with turning proportions. Both studies were able to describe kinematic waves arising from a network intersection but also subject to significant shortcomings: (i) All links are assumed to have the same fundamental diagram in both studies; (ii) In \citep{holden1995unidirection}, the entropy conditions are pragmatic and lack of physical interpretations; and (iii) In \citep{coclite2005network}, the Riemann problem can only be uniquely solved for junctions with no fewer downstream links, and hence their results do not apply to a merging junction.

In this paper, we are interested in studying continuous kinematic wave models of merging traffic flow which are consistent with discrete CTM with various distribution schemes. Here we consider a merge network with $m\geq 2$ upstream links and one downstream link, as shown in \reff{merge20080729figure.22}. In this network, there are $m+1$ links, $m$ origin-destination pairs, and $m$ paths. 
In the continuous kinematic wave models of this network, traffic dynamics on each link are described by the LWR model \refe{lwr}, and, in addition, an entropy condition based on various distribution schemes is used to pick out unique physical solutions \citep{ansorge1990entropy}. Note that traffic dynamics on the whole network cannot be modeled by either one-dimensional or two-dimensional conservation laws in closed forms, since vehicles of different commodities interact with each other on the downstream link. Similar to the LWR model, it is not possible to obtain analytical solutions for general initial and boundary conditions, and we usually resort to analyzing continuous kinematic wave solutions to Riemann problems. Studies on the Riemann problem are helpful for understanding fundamental characteristics of the corresponding merge model and constructing numerical solutions. In the Riemann problem, all links are homogeneous and infinitely long; for link $a=1,\cdots, m+1$, we assume that its flow-density relation is $q_a=Q_a(\r_a)$, critical density $\r_{c,a}$, and its capacity $C_a$; and upstream link $i=1,\cdots,m$ and downstream link $m+1$ have constant initial conditions:
\bqn
\r_i(x_i,0)&=&\r_{i}, \: x_i\in(-\infty,0), \quad i=1,\cdots,m \label{linkini1} \\
\r_{m+1}(x_{{m+1}},0)&=&\r_{{m+1}}, \: x_{{m+1}}\in(0,+\infty). \label{linkini2}
\eqn
To pick out physical kinematic wave solutions to the Riemann problem, here we use various distribution schemes of CTM as our entropy conditions, which are different from the optimization-based entropy condition in \citep{holden1995unidirection}. In this sense, the resulted models can be considered as continuous versions of CTM with various distribution schemes. 

\bfg\bc
\includegraphics[height=4cm]{merge20080729figure.22}\caption{An illustration of a merge network}\label{merge20080729figure.22}
\ec\efg

In our study, the upstream links can be mainline freeways or on-ramps with different characteristics, and our solutions are physically meaningful and consistent with discrete CTM of merging traffic. In addition, we follow a new framework developed in \citep{jin2009sd} for solving the Riemann problem for inhomogeneous LWR model at a linear junction. In this framework, we first construct Riemann solutions of a merge model in the supply-demand space: on each branch, a stationary state will prevail after a long time, and an interior state may exist but does not take any space in the continuous solution and only shows up in one cell in the numerical solutions \citep{vanleer1984upwind}. After deriving admissible solutions for upstream and downstream stationary and interior states, we introduce an entropy condition consistent with various discrete merge models. We then prove that stationary states and boundary fluxes are unique for given upstream demands and downstream supply, and interior states may not be unique but inconsequential. 
Then, the kinematic wave on a link is uniquely determined by the corresponding LWR model with the stationary state and the initial state as initial states.

In this study, we do not attempt to develop new simulation models. Rather, we attempt to understand kinematic wave solutions of continuous merge models consistent with various CTM. The theoretical approach is not meant to replace numerical methods in CTM but to provide insights on the formation and dissipation of shock and rarefaction waves for the Riemann problem. In this study, we choose to solve the Riemann problem in supply-demand space due to its mathematical tractability. In fact, the benefits of working in the supply-demand space have been demonstrated by \citet{daganzo1995ctm} and \citet{lebacque1996godunov} when they developed discrete supply-demand models of network vehicular traffic, i.e., CTM. In their studies, the introduction of the concepts of demand and supply enabled a simple and intuitive formulation of merge and diverge models. But different from these studies, which focused on discrete simulation models, this study further introduces supply-demand diagram and discusses analytical solutions of continuous merge models.

The rest of the paper is organized as follows. In Section 2, we introduce an analytical framework for solving the kinematic waves of the Riemann problem with jump initial conditions in supply-demand space. In particular, we derive traffic conservation conditions, admissible conditions of stationary and initial states, and additional entropy conditions consistent with discrete merge models. In Section 3, we solve stationary states and boundary fluxes for both fair and constant merge models. In Section 4, we introduce two invariant merge models. In Section 5, we demonstrate the validity of the proposed analytical framework with numerical examples. In Section 6, we summarize our findings and present some discussions.

\section{{An analytical framework}}
For link $a=1,\cdots, m+1$, we define the following demand and supply functions with all subscript $a$ suppressed \citep{engquist1980calculation,daganzo1995ctm,lebacque1996godunov}
\bqn
D(\r)&=&Q(\min\{\r,\r_c\})=\cas{{ll}Q(\r),&\m{if } \r\leq \r_c\\C,&\m{if }\r\geq \r_c},\nonumber \\&=&\int_0^\r \chi(s) Q'(s) ds=\int_0^\r \max\{Q'(s),0\} ds \label{def:d}\\
S(\r)&=&Q(\max\{\r,\r_c\})=\cas{{ll}Q(\r),&\m{if } \r\geq \r_c\\C,&\m{if }\r\leq \r_c},\nonumber\\&=&C+\int_0^\r (1- \chi(s)) Q'(s) ds=C+\int_0^\r \min\{Q'(s),0\} ds, \label{def:s}
\eqn
where $\chi(\r)$ equals 1 iff $Q'(\r)\geq 0$ and equals 0 otherwise. 

Unlike many existing studies, in which traffic states are considered in $\r$-$q$ space, we represent a traffic state in supply-demand space as $U=(D,S)$. 
For the demand and supply functions in \refe{def:d} and \refe{def:s}, we can see that $D$ is non-decreasing with $\r$ and $S$ non-increasing.
Thus $D\leq C$, $S\leq C$, $\max \{D,S\}=C$, and flow-rate $q(U)=\min\{D,S\}$. In addition, $D=S=C$ iff traffic is critical; $D<S=C$ iff traffic is strictly under-critical (SUC); $S<D=C$ iff traffic is strictly over-critical (SOC). Therefore, state $U=(D,S)$ is under-critical (UC), iff $S=C$, or equivalently $D\leq S$; State $U=(D,S)$ is over-critical (OC), iff $D=C$, or equivalently $S\leq D$.

In \reff{fig:fd-ds}(b), we draw a supply-demand diagram for the two fundamental diagrams in \reff{fig:fd-ds}(a).
 On the dashed branch of the supply-demand diagram, traffic is UC and $U=(D,C)$ with $D\leq C$; on the solid branch, traffic is OC and $U=(C,S)$ with $S\leq C$. Compared with the fundamental diagram of a road section, the supply-demand diagram only considers its capacity $C$ and criticality, but not other detailed characteristics such as  critical density, jam density, or the shape of the fundamental diagram.
That is, different fundamental diagrams can have the same demand-supply diagram, as long as they have the same capacity and are unimodal. 
However, given a demand-supply diagram and its corresponding fundamental diagram, the points are one-to-one mapped. 
  
 \bfg\bc $\ba{c@{\hspace{0.3in}}c}
\includegraphics[height=1.8in]{sta20061026figure.2} &
\includegraphics[height=1.8in]{sta20061026figure.1} \\
\multicolumn{1}{c}{\mbox{\bf (a)}} &
    \multicolumn{1}{c}{\mbox{\bf (b)}}
\ea$ \ec \caption{Fundamental diagrams and their corresponding supply-demand diagrams }\label{fig:fd-ds} \efg
 
In supply-demand space, initial conditions in \refe{linkini1} and \refe{linkini2} are equivalent to \footnote{In this section, $i=1,\cdots,m$ if not otherwise mentioned.}
\bqn
U_i(x_i,0)&=&(D_i,S_i), \quad x_i\in(-\infty, 0), \label{dslinkini1}  \\
U_{m+1}(x_{m+1},0)&=&(D_{m+1},S_{m+1}), \quad x_{m+1}\in(0, +\infty).\label{dslinkini2}
\eqn
In the solutions of the Riemann problem with initial conditions (\ref{dslinkini1}-\ref{dslinkini2}), a shock wave or a rarefaction wave could initiate on a link from the merging junction $x=0$, and traffic states on all links become stationary after a long time. That is, in solutions to the Riemann problem, stationary states prevail all links after a long time. At the boundary, there can also exist interior states \citep{vanleer1984upwind,bultelle1998shock}, which take infinitesimal space \footnote{In numerical solutions, the interior states exist in one cell.}. We denote the stationary states on upstream link $i$ and  downstream link ${m+1}$ by $U_i^-$ and $U_{m+1}^+$, respectively. We denote the interior states on links $i$ and ${m+1}$ by $U_i(0^-,t)$ and $U_{m+1}(0^+,t)$, respectively. The structure of Riemann solutions on upstream and downstream links are shown in \reff{fig:riemannstructure}, where arrows illustrate the directions of possible kinematic waves. 
Then the kinematic wave on upstream link $i$ is the solution of the corresponding LWR model with initial left and right conditions of $U_i$ and $U_i^-$, respectively. Similarly, the kinematic wave on downstream link ${m+1}$ is the solution of the corresponding LWR model with initial left and right conditions of $U_{m+1}^+$ and $U_{m+1}$, respectively. 

\bfg\bc $\ba{c@{\hspace{0.3in}}c}
\includegraphics[width=2.9in]{merge20080729figure.4} &
\includegraphics[width=2.9in]{merge20080729figure.5} \\
\multicolumn{1}{c}{\mbox{\bf (a)}} &
    \multicolumn{1}{c}{\mbox{\bf (b)}}
\ea$ \ec \caption{The structure of Riemann solutions: (a) Upstream link $i$; (b) Downstream link $m+1$}\label{fig:riemannstructure} \efg

We denote $q_{i\to m+1}$ as the flux from link $i$ to link $m+1$ for $t>0$. The fluxes are determined by the stationary states: the out-flux of link $i$ is $q_i=q(U_i^-)$, and the in-flux of link $m+1$ is $q_{m+1}=q(U_{m+1}^+)$. Furthermore, from traffic conservation at a merging junction, we have 
\bqn
q_{i\to m+1}&=&q_i=q(U_i^-), \quad q_{m+1}=q(U_{m+1}^+)=\sum_{i=1}^{m} q(U_i^-). \label{trafficconservation}
\eqn

\subsection{Admissible stationary and interior states}

As observed in \citep{holden1995unidirection,coclite2005network}, the speed of a kinematic wave on an upstream link cannot be positive, and that on a downstream link cannot be negative. We have the following admissible conditions on stationary states.
\begin{theorem} [Admissible stationary states] \label{admissibless} For initial conditions in \refe{dslinkini1} and \refe{dslinkini2},  stationary states are admissible if and only if 
\bqn
U_i^-&=&(D_i,C_i) \m{ or } (C_i,  S_i^-), \label{upstreamss}
\eqn
where $S_i^-<D_i$ , and 
\bqn
U_{m+1}^+&=&(C_{m+1},S_{m+1}) \m{ or } (D_{m+1}^+, C_{m+1}), \label{downstreamss}
\eqn
where $D_{m+1}^+<S_{m+1}$ .
\end{theorem}
The proof is quite straightforward and omitted here. The regions of admissible upstream stationary states in both supply-demand and fundamental diagrams are shown in \reff{fig:ssupstream}, and the regions of admissible downstream stationary states are shown in \reff{fig:ssdownstream}. In \reff{fig:ssupstream}(a) and (b), the initial upstream condition is SUC with $D_i<C_i=S_i$, and feasible stationary states are $(D_i,C_i)$, when there is no wave, or $(C_i,S_i^-)$ with $S_i^-<D_i$, when a back-traveling shock wave emerges on the upstream link $i$; In \reff{fig:ssupstream}(c) and (d), the initial upstream condition is OC with $S_i\leq D_i=C_i$, and any OC stationary state is feasible, when a back-traveling shock or rarefaction wave emerges on the upstream link $i$. In contrast, in \reff{fig:ssdownstream}(a) and (b), the initial downstream condition is UC with $D_{m+1}\leq S_{m+1}=C_{m+1}$, and any UC stationary state is feasible, when a forward-traveling shock or rarefaction wave emerges on the downstream link $m+1$; In \reff{fig:ssdownstream}(c) and (d), the initial downstream condition is SOC with $S_{m+1}<D_{m+1}=C_{m+1}$, and feasible stationary states are $(C_{m+1}, S_{m+1})$, when there is no wave, or $(D_{m+1}^+, C_{m+1})$ with $D_{m+1}^+<S_{m+1}$, when a forward-traveling shock wave emerges on the downstream link $m+1$. Here the types of kinematic waves and the signs of the wave speeds can be determined in the supply-demand diagram, but the absolute values of the wave speeds have to be determined in the fundamental diagram.

\bfg\bc $\ba{c@{\hspace{0.3in}}c}
\includegraphics[height=1.8in]{merge20080729figure.6} &
\includegraphics[height=1.8in]{merge20080729figure.7} \\
\multicolumn{1}{c}{\mbox{\bf (a)}} &
    \multicolumn{1}{c}{\mbox{\bf (b)}}\\
\includegraphics[height=1.8in]{merge20080729figure.8} &
\includegraphics[height=1.8in]{merge20080729figure.9} \\
\multicolumn{1}{c}{\mbox{\bf (c)}} &
    \multicolumn{1}{c}{\mbox{\bf (d)}}	
\ea$ \ec \caption{Admissible stationary states for upstream link $i$: marked by black dots}\label{fig:ssupstream} \efg

\bfg\bc $\ba{c@{\hspace{0.3in}}c}
\includegraphics[height=1.8in]{merge20080729figure.10} &
\includegraphics[height=1.8in]{merge20080729figure.11} \\
\multicolumn{1}{c}{\mbox{\bf (a)}} &
    \multicolumn{1}{c}{\mbox{\bf (b)}}\\
\includegraphics[height=1.8in]{merge20080729figure.12} &
\includegraphics[height=1.8in]{merge20080729figure.13} \\
\multicolumn{1}{c}{\mbox{\bf (c)}} &
    \multicolumn{1}{c}{\mbox{\bf (d)}}	
\ea$ \ec \caption{Admissible stationary states for downstream link $m+1$: marked by black dots}\label{fig:ssdownstream} \efg

{\em Remark 1.} $U_i^-=U_i$ and $U_{m+1}^+=U_{m+1}$ are always admissible. In this case, the stationary states are the same as the corresponding initial states, and there are no waves.

{\em Remark 2.}
Out-flux $q_i=\min\{D_i^-,S_i^-\}\leq D_i$ and in-flux $q_{m+1}=\min\{D_{m+1}^+,S_{m+1}^+\}\leq S_{m+1}$. That is, $D_i$ is the maximum sending flow and $S_{m+1}$ is the maximum receiving flow in the sense of \citep{daganzo1994ctm,daganzo1995ctm}.

{\em Remark 3.} The ``invariance principle" in \citep{lebacque2005network} can be interpreted as follows: if $D_i^-=C_i$, then $q(U_i^-)<D_i$; if $S_{m+1}^+=C_{m+1}$, then $q(U_{m+1}^+)<S_i$. We can see that Theorem \ref{admissibless} is consistent with the ``invariance principle".

\begin{corollary} \label{flux2stationary} For the upstream link $i$, $q_i\leq D_i$; $q_i<D_i$ if and only if $U_i^-=(C_i,q_i)$, and $q_i=D_i$ if and only if $U_i^-=(D_i, C_i)$. For the downstream link $m+1$, $q_{m+1}\leq S_{m+1}$; $q_{m+1}<S_{m+1}$ if and only if $U_{m+1}^+=(q_{m+1},C_{m+1})$, and $q_{m+1}=S_{m+1}$ if and only if $U_{m+1}^+=(C_{m+1},S_{m+1})$. That is, given out-fluxes and in-fluxes, the stationary states can be uniquely determined.
\end{corollary}

For interior states, the waves of the Riemann problem on link $i$ with left and right initial conditions of $U_i^-$ and  $U_i(0^-,t)$ cannot have negative speeds. Similarly, the waves of the Riemann problem on link $m+1$ with left and right initial conditions of $U_{m+1}(0^+,t)$ and $U_{m+1}^+$ cannot have positive speeds. Therefore, interior states $U_i(0^-,t)$ and $U_{m+1}(0^+,t)$ should satisfy the following admissible conditions. 
\begin{theorem}[Admissible interior states]\label{intersta}
For stationary states $U_i^-$ and $U_{m+1}^+$, interior states $U_i(0^-,t)$ and $U_{m+1}(0^+,t)$ are admissible if and only if  
\bqn
U_i(0^-,t)&=&\cas{{ll}(C_i,S_i^-)=U_i^-, &\m{when }U_i^-\m{ is SOC; i.e., }S_i^-< D_i^-=C_i \\(D_i(0^-,t),  S_i(0^-,t)), & \m{when }U_i^-\m{ is UC; i.e., }D_i^-\leq S_i^-=C_i} \label{upstreamis}
\eqn
where $S_i(0^-,t) \geq D_i^-$ , and 
\bqn
U_{m+1}(0^+,t)&=&\cas{{ll}(D_{m+1}^+,C_{m+1})=U_{m+1}^+,&\m{when }U_{m+1}^+\m{ is SUC; i.e., }D_{m+1}^+<S_{m+1}^+=C_{m+1} \\ (D_{m+1}(0^+,t), S_{m+1}(0^+,t)), &\m{when }U_{m+1}^+\m{ is OC; i.e., }S_{m+1}^+\leq D_{m+1}^+=C_{m+1} } \label{downstreamis}
\eqn
where $D_{m+1}(0^+,t)\geq S_{m+1}^+$ .
\end{theorem}
The proof is quite straightforward and omitted here. The regions of admissible upstream interior states in both supply-demand and fundamental diagrams are shown in \reff{fig:isupstream}, and the regions of admissible downstream interior states are shown in \reff{fig:isdownstream}. From the figures, we can also determine the types and traveling directions of waves with given stationary and interior states on all links, but these waves are suppressed and cannot be observed. Note that, however, we are able to observe possible interior states in numerical solutions.

\bfg\bc $\ba{c@{\hspace{0.3in}}c}
\includegraphics[height=1.8in]{merge20080729figure.14} &
\includegraphics[height=1.8in]{merge20080729figure.15} \\
\multicolumn{1}{c}{\mbox{\bf (a)}} &
    \multicolumn{1}{c}{\mbox{\bf (b)}}\\
\includegraphics[height=1.8in]{merge20080729figure.16} &
\includegraphics[height=1.8in]{merge20080729figure.17} \\
\multicolumn{1}{c}{\mbox{\bf (c)}} &
    \multicolumn{1}{c}{\mbox{\bf (d)}}	
\ea$ \ec \caption{Admissible interior states for upstream link $i$: marked by black dots}\label{fig:isupstream} \efg

\bfg\bc $\ba{c@{\hspace{0.3in}}c}
\includegraphics[height=1.8in]{merge20080729figure.18} &
\includegraphics[height=1.8in]{merge20080729figure.19} \\
\multicolumn{1}{c}{\mbox{\bf (a)}} &
    \multicolumn{1}{c}{\mbox{\bf (b)}}\\
\includegraphics[height=1.8in]{merge20080729figure.20} &
\includegraphics[height=1.8in]{merge20080729figure.21} \\
\multicolumn{1}{c}{\mbox{\bf (c)}} &
    \multicolumn{1}{c}{\mbox{\bf (d)}}	
\ea$ \ec \caption{Admissible interior states for downstream link $m+1$: marked by black dots}\label{fig:isdownstream} \efg

{\em Remark 1.} Note that $U_i(0^-,t)=U_i^-$ and $U_{m+1}(0^+,t)=U_{m+1}^+$ are always admissible. In this case, the interior states are the same as the stationary states.

\subsection{Entropy conditions consistent with discrete merge models}
In order to uniquely determine the solutions of stationary states, we introduce a so-called entropy condition in interior states as follows:
\bqn
q_i&=&F_i(U_1(0^-,t),\cdots, U_m(0^-,t), U_{m+1}(0^+,t)).  \label{sdentropy}
\eqn
The entropy condition can also be written as
\bqs
q_{m+1}&=&F(U_1(0^-,t),\cdots, U_m(0^-,t), U_{m+1}(0^+,t)),\\
q_i&=& \alpha_i q_{m+1},
\eqs
where $\alpha_i \in[0, 1]$ and $\sum_{i=1}^m \alpha_i=1$.
We can see that the entropy condition uses ``local" information in the sense that it determines boundary fluxes from interior states. In the discrete merge model, the entropy condition is used to determine boundary fluxes from cells contingent to the merging junction.
Thus, $F_i(U_1(0^-,t),\cdots, U_m(0^-,t), U_{m+1}(0^+,t))$ in \refe{sdentropy} can be considered as local, discrete flux functions.

Here the boundary fluxes can be obtained with existing merge models.
In \citep{daganzo1995ctm}, $F(U_1(0^-,t),\cdots, U_m(0^-,t), U_{m+1}(0^+,t))$ was proposed to solve the following local optimization problem
\bqn
\max_{U_i^-,U_{m+1}^+,U_i(0^-,t),U_{m+1}(0^+,t)} \{q_{m+1} \} \label{optimizationentropy}
\eqn
subject to 
\bqs
&&q_i \leq  D_i(0^-,t),\\
&&q_{m+1} \leq  S_{m+1}(0^+,t).
\eqs
Thus, we obtain the total flux as
\bqs
F(U_1(0^-,t),\cdots, U_m(0^-,t), U_{m+1}(0^+,t))=\min\{\sum_{i=1}^m D_i(0^-,t), S_{m+1}(0^+,t) \}.
\eqs
The proportions, $\alpha_i$, can be determined by so-called ``distribution schemes", which distribute the total flux to each upstream link.
In \citep{newell1993sim}, an on-ramp is given total priority and can send its maximum flow. In \citep{daganzo1995ctm}, a priority-based scheme was proposed. In \citep{lebacque1996godunov}, a general scheme was proposed, and it was suggested to distribute out-fluxes according to the number of lanes of upstream links. In \citep{banks2000metering}, an on-ramp is given total priority, but its flow is also restricted by the metering rate. In \citep{jin2003merge}, a fair scheme is proposed to distribute out-fluxes according to local traffic demands of upstream links. In \citep{ni2005merge}, a fair share of the downstream supply is assigned to each upstream proportional to its capacity, and  out-fluxes are then determined by comparing the corresponding fair shares and demands.
With these distribution schemes, the entropy condition, \refe{sdentropy}, yields unique solutions of boundary fluxes with given interior states.

\commentout{
In the literature, 

Therefore, we can introduce the following entropy condition for the Riemann problem with initial conditions in (\ref{dslinkini1}-\ref{dslinkini2}):
\bi
\item [] The boundary fluxes through a merging junction are the solutions of the optimization problem in \refe{sdentropy}.
\ei
}

\subsection{Summary of the solution framework}
To solve the Riemann problem with the initial conditions in \refe{dslinkini1}-\refe{dslinkini2}, we will first find stationary and interior states that satisfy the aforementioned entropy condition, admissible conditions, and traffic conservation equations. Then the kinematic wave on each link will be determined by the Riemann problem of the corresponding LWR model with initial and stationary states as initial conditions \citep{jin2009sd}. 
Here we will only focus on solving the stationary states on all links.
We can see that the feasible domains of stationary and interior states are independent of the upstream supply, $S_i$, and the downstream demand, $D_{m+1}$. That is, the same upstream demand and downstream supply will yield the same solutions of stationary and interior states. However, the upstream and downstream wave types and speeds on each can be related to $S_i$ as shown in \reff{fig:ssupstream}(d) and $D_{m+1}$ as shown in \reff{fig:ssdownstream}(d).

\section{{Solutions of two merge models}}
In this paper, we solve the Riemann problem for a merging junctions with two upstream links; i.e., $m=2$. 
Different merge models have different entropy conditions \refe{sdentropy}. In this section, we consider two entropy conditions, i.e., two merge models.
Here we attempt to find the relationships between the boundary fluxes and the initial conditions.
\bqn
q_i=\hat F_i(U_1,U_2, U_3).
\eqn
In contrast to local, discrete flux functions $F_i(U_1(0^-,t),U_2(0^-,t), U_3(0^+,t))$, $\hat F_i(U_1,U_2, U_3)$ can be considered as global, continuous.
With the global, continuous fluxes, we can find stationary states from Corollary \ref{flux2stationary}. With the solution framework in the preceding section, we can then find the kinematic waves of the Riemann problem with initial conditions $(U_1, U_2, U_3)$. 

\subsection{A fair merge model}
We consider the fair merging rule proposed in \citep{jin2003merge}, in which
\bqn
\a_i=\frac{D_i(0^-,t)}{D_1(0^-,t)+D_2(0^-,t)},\quad i=1,2. \label{fairdistribution}
\eqn
Obviously, the optimization problem of \refe{sdentropy} subject to the fair merging rule yields the following solutions
\bqn
q_i&=&\min\{1,\frac{S_3(0^+,t)}{D_1(0^-,t)+D_2(0^-,t)}\} D_i(0^-,t),\quad i=1,2. \label{sdentropy:fair}
\eqn
Thus in the Riemann solutions, stationary and interior states have to satisfy \refe{sdentropy:fair}, traffic conservation, and the corresponding admissible conditions.

\begin{theorem} \label{thm:fairmerg2} For the Riemann problem with the initial conditions in \refe{dslinkini1} and \refe{dslinkini2},  stationary and interior states satisfying the entropy condition in \refe{sdentropy:fair}, traffic conservation equations, and the corresponding admissible conditions are the following:
\ben
\item When $D_1+D_2 < S_{3}$, $U_i^-=U_i(0^-,t)=(D_i,C_i)$ ($i=1,2$) and $U_{3}^+=U_3(0^+,t)=(D_1+D_2, C_{3})$;
\item When $D_1+D_2 = S_{3}$, $U_i^-=U_i(0^-,t)=(D_i,C_i)$ ($i=1,2$), $U_{3}^+=(C_{3}, S_{3})$, $U_3(0^+,t)=(C_{3}, S_{3})$ or $(D_3(0^+,t), S_3(0^+,t))$ with $D_3(0^+,t)\geq S_3$ and $S_3(0^+,t)>S_3$ when $S_3<C_3$;
\item When $D_i> \frac{C_i}{C_1+C_2}S_3$ ($i=1,2$), $U_{3}^+=U_3(0^+,t)=(C_{3}, S_{3})$, and $U^-_i=U_i(0^-,t)=(C_i, \frac {C_i}{C_1+C_2} S_3 )$.
\item When $D_1+D_2 > S_{3}$ and $D_i \leq \frac{C_i}{C_1+C_2}S_3$ ($i,j=1$ or 2 and $i\neq j$), $U_{3}^+=U_3(0^+,t)=(C_{3}, S_{3})$, $U^-_i=(D_i, C_i)$, $U_i(0^-,t)=(\frac{C_j}{S_3-D_i} D_i,C_i)$, and $U^-_j=U_j(0^-,t)=(C_j, S_3-D_i)$.
\een
\end{theorem}
The proof of the theorem is given in Appendix A. 
In \reff{fig:merg2}, we demonstrate how to obtain the pair of $(D_1,D_2)$ from initial states $U_1$ and $U_2$ and how to obtain the line of $S_3-S_3$ in the supply-demand diagrams. According to the relationship between $(D_1,D_2)$ and $S_3$, the domain of $(D_1,D_2)$ can be divided into four regions. In \reff{fig:type1}, we further demonstrate how to obtain stationary states for all these four regions. The results are consistent with those in \citep{daganzo1996gridlock} in principle: for initial conditions in region I, there are no queues on upstream links; for initial states in region II, there is a queue on link 2; for initial states in region III, there are queues on both links; and, for initial states in region IV, there is a queue on link 1.

\bfg
\bc\includegraphics{sta20061026figure.12}\ec
\caption{Four types of solutions for a fair merge model}\label{fig:merg2}
\efg

\bfg \bc $\ba{c@{\hspace{0.3in}}c}
\includegraphics[height=2in]{sta20061026figure.14} &
\includegraphics[height=2in]{sta20061026figure.15}\\
\multicolumn{1}{c}{\mbox{\bf (a)}} &
    \multicolumn{1}{c}{\mbox{\bf (b)}} \\
    \includegraphics[height=2in]{sta20061026figure.16} &
\includegraphics[height=2in]{sta20061026figure.17}\\
\multicolumn{1}{c}{\mbox{\bf (c)}} &
    \multicolumn{1}{c}{\mbox{\bf (d)}} 
\ea$
\ec
\caption{Solutions of stationary states for different initial conditions for the fair merge model}\label{fig:type1} \efg

\begin{corollary} \label{fairflux}
For the Riemann problem with initial conditions in \refe{dslinkini1} and \refe{dslinkini2}, the boundary fluxes satisfying the entropy condition in \refe{sdentropy:fair}, traffic conservation equations, and the corresponding admissible conditions are in the following:
\ben
\item When $D_1+D_2 \leq S_{3}$, $q_i=D_i$ ($i=1,2$) and $q_3=D_1+D_2$;
\item When $D_i > \frac{C_i}{C_1+C_2}S_3$ ($i=1,2$), $q_i=\frac{C_i}{C_1+C_2}S_3$ and $q_3=S_3$;
\item When $D_1+D_2 > S_{3}$ and $D_i \leq \frac{C_i}{C_1+C_2}S_3$ ($i,j=1$ or 2 and $i\neq j$),  $q_i=D_i$, $q_j=S_3-D_i$, and $q_3=S_3$.
\een
That is, for $i,j=1$ or 2 and $i\neq j$,
\bqn
q_i&=&\min\{D_i, \max\{ S_3-D_j,\frac{S_3}{C_1+C_2} C_i\}\}. \label{sdentropy:fair2}
\eqn
\end{corollary}
The solutions of fluxes in four different regions are shown in \reff{merge20080729figure.1}, in which the dots represent the initial conditions in $(D_1,D_2)$, and the end points of the arrows represent the solutions of fluxes $(q_1,q_2)$. 
The fluxes in \refe{sdentropy:fair2} are exactly the same as in \citep{ni2005merge}. 
In this sense, the discrete merge model with \refe{sdentropy:fair} converges to the merge model in \citep{ni2005merge}. 

\bfg
\bc\includegraphics[height=8cm]{merge20080729figure.1}\ec
\caption{Solutions of fluxes for the fair merge model}\label{merge20080729figure.1}
\efg

\begin{corollary}
If $U_i$ ($i=1,2$) and $U_3$ satisfy
\bqs
\min\{D_i,S_i\}&=&\min\{D_i, \max\{ S_3-D_j,\frac{S_3}{C_1+C_2} C_i\}\},\\
\min\{D_3,S_3\}&=&\min\{D_1, \max\{ S_3-D_2,\frac{S_3}{C_1+C_2} C_1\}\}+\min\{D_2, \max\{ S_3-D_1,\frac{S_3}{C_1+C_2} C_2\}\}\\&=&\min\{D_1+D_2,S_3\},
\eqs
then the unique stationary states are the same as the initial states, and traffic dynamics at the merging junction are stationary.
\end{corollary}
{\em Proof}. The results follow from Theorem \ref{thm:fairmerg2}, from which we can see that stationary states $U_i^-$ and $S_3^+$ satisfy
\bqs
\min\{D_i^-,S_i^-\}&=&\min\{D_i^-, \max\{ S_3^+-D_j^-,\frac{S_3^+}{C_1+C_2} C_i\}\},\\
\min\{D_3^+,S_3^+\}&=&\min\{D_1^-, \max\{ S_3^+-D_2^-,\frac{S_3^+}{C_1+C_2} C_1\}\}+\min\{D_2^-, \max\{ S_3^+-D_1^-,\frac{S_3^+}{C_1+C_2} C_2\}\}\\&=&\min\{D_1^-+D_2^-,S_3^+\}.
\eqs
\eop

\subsection{A constant merge model}
We consider a merge model proposed in \citep{lebacque1996godunov}, in which
\bqn
q_i=\min\{D_i(0^-,t), \alpha_i S_{3}(0^+,t)\}, \quad  i=1,2 \label{sdentropy:constant}
\eqn
where $\alpha_i$ are constant distribution proportions, $\alpha_i\in [0,1]$, and $\sum_i \alpha_i=1$.

\begin{theorem} \label{thm:constmerg2} For the Riemann problem with the initial conditions in \refe{dslinkini1} and \refe{dslinkini2}, boundary fluxes satisfying the entropy condition in \refe{sdentropy:constant}, traffic conservation equations, and the corresponding admissible conditions are the following:
\ben
\item When $D_1+D_2 < S_{3}$ and $D_i\leq \a_i C_3$ ($i=1,2$), $q_i=D_i$ and $q_3=D_1+D_2$;
\item When $D_i >\a_i C_3$ and $D_j<S_3-\a_i C_3$ ($i,j=$ 1 or 2 and $i\neq j$), $q_i=\a_i C_3$, $q_j=D_j$, and $q_3=\a_i C_3+D_j$.
\item When $D_1+D_2\geq S_3$, $S_3-\a_j C_3 \leq D_i \leq \a_i S_3$ ($i,j=$1 or 2 and $i\neq j$), $q_i=D_i$, $q_j=S_3-D_i$, and $q_3=S_3$. 
\item When $D_i \geq \alpha_i S_3$ ($i=1,2$), $q_i=\a_i S_3$, and $q_3=S_3$. 
\een
That is, for $i,j=1$ or 2 and $i\neq j$,
\bqn
q_i&=&\min\{D_i, \a_i C_3, \max\{ S_3-D_j, \a_i S_3 \}\}. \label{sdentropy:constant2}
\eqn
\end{theorem}
The proof of the theorem is in Appendix B. The solutions of fluxes in six different regions are shown in \reff{merge20080729figure.3}, in which the starting point of an arrow represents initial conditions, $(D_1,D_2)$, and the ending point represents the solution of fluxes, $(q_1,q_2)$. With the fluxes, we can easily find all stationary and interior states as in Theorem \ref{thm:fairmerg2}. 
Compared with the fair merging rule, the constant merging rule yield sub-optimal solutions in regions II and VI in \reff{merge20080729figure.3}, in which $q_3<\min\{D_1+D_2,S_3\}$. That is, this merge model does not satisfy the optimization entropy condition in \refe{optimizationentropy}.

\bfg
\bc\includegraphics[height=8cm]{merge20080729figure.3}\ec
\caption{Solutions of fluxes for the constant merge model}\label{merge20080729figure.3}
\efg

\commentout{The equivalence between \refe{sdentropy:constant2} and \refe{sdentropy:fair2}, when $\a_i=\frac{D_i}{D_i+D_j}$? Probably not, since the proof is demand independent. Stationary traffic dynamics. Other solutions?}

\section{Invariant merge models}
For the two merge models in the preceding section, the local, discrete flux functions $F_i(\cdot, \cdot, \cdot)$ are not the same as the global, continuous flux functions $\hat F_i(\cdot, \cdot, \cdot)$. 
That is, boundary fluxes computed from local flux functions are not constant. 
In this sense, the fair and constant merge models are not invariant. In the following, we propose two invariant merge models, in which $F_i(\cdot, \cdot, \cdot)$ and $\hat F_i(\cdot, \cdot, \cdot)$ have the same functional form.

\subsection{An optimal invariant merge model}
We consider the following priority-based merge model \citep{daganzo1995ctm,daganzo1996gridlock}, in which $i,j=1, \m{ or }2$ and $i\neq j$
\bqn
q_i=\min\{D_i(0^-,t), \max\{S_3(0^+,t)-D_j(0^-,t), \alpha_i S_3(0^+,t)\}\}, \quad \label{sdentropy:priority}
\eqn
where $\alpha_i$ are priority distribution proportions, $\alpha_i\in[0,1]$, and $\sum_i \alpha_i=1$.

\begin{theorem} \label{thm:invariant1}
For the Riemann problem with the initial conditions in \refe{dslinkini1} and \refe{dslinkini2}, boundary fluxes satisfying the entropy condition in \refe{sdentropy:priority}, the traffic conservation equations, and the corresponding admissible conditions are the following:
\bqn
q_i=\min\{D_i, \max\{S_3-D_j, \alpha_i S_3\}\}, \quad i,j=1, \m{ or }2,\m{ and } i\neq j \label{sdentropy:priority2}
\eqn
\end{theorem}

The proof of the theorem is given in Appendix C. Since $F_i(\cdot,\cdot,\cdot)=\hat F_i(\cdot,\cdot,\cdot)$, the merge model \refe{sdentropy:priority} is invariant.
From the proof, we can see that the merging rule is optimal in the sense that $q_3=\min\{D_1+D_2,S_3\}$. The solutions of fluxes in four different regions are shown in \reff{merge20080729figure.2}, in which the starting point of an arrow represents the initial demands, $(D_1,D_2)$, and the end point represents the solutions of fluxes, $(q_1,q_2)$. 
Note that, when $\a_i=\frac{C_i}{C_1+C_2}$ ($i=1,2$), then the distribution scheme is the same as \refe{sdentropy:fair2}. 
Therefore, the merge model in \citep{ni2005merge} is invariant.

\bfg
\bc\includegraphics[height=8cm]{merge20080729figure.2}\ec
\caption{Solutions of fluxes for the priority-based merge model}\label{merge20080729figure.2}
\efg

When $\a_i>\frac{C_i}{C_1+C_2}$, the merge model \refe{sdentropy:priority} gives higher priority to vehicles from upstream link $i$. 
Thus, \refe{sdentropy:priority} is another representation of the priority merge model proposed in \citep{daganzo1995ctm,daganzo1996gridlock}.
For an extreme merge model with $\a_1=1$ and $\a_2=0$ \citep{newell1993sim,banks2000metering}, vehicles from link 1 is given absolute priorities. In this case, we have
\bqs
q_1&=&\min\{D_1, S_3\},\\
q_2&=&\min\{D_2, \max\{S_3-D_1, 0\}\}.
\eqs

\subsection{A suboptimal invariant merge model}
Similarly, the continuous version of the constant merge model will lead to another invariant merge model, in which
\bqn
q_i=\min\{D_i(0^-,t), \a_i C_3, \max\{S_3(0^+,t)-D_j(0^-,t), \alpha_i S_3(0^+,t)\}\}, \quad i,j=1, \m{ or }2,\m{ and } i\neq j \label{invariant:constant}
\eqn
with $\alpha_i\in[0,1]$, and $\sum_i \alpha_i=1$.

\begin{theorem} \label{thm:invariant2}
For the Riemann problem with initial conditions in \refe{dslinkini1} and \refe{dslinkini2}, the boundary fluxes satisfying the entropy condition in \refe{sdentropy:priority}, traffic conservation equations, and the corresponding admissible conditions are in the following:
\bqn
q_i=\min\{D_i, \a_i C_3, \max\{S_3-D_j, \alpha_i S_3\}\}, \quad i,j=1, \m{ or }2,\m{ and } i\neq j \label{invariant:constant2}
\eqn
\end{theorem}

The proof of the theorem is omitted. The merge model \refe{invariant:constant} is obviously invariant. In addition, the solutions of fluxes in six different regions are the same as in \reff{merge20080729figure.3}. We can see that this merge model is suboptimal.
Assuming that both links 1 and 3 are the mainline freeway with capacity $C_1=C_3$ and $D_1=S_3=C_1$, and link 2 is a metered on-ramp with a metering rate $D_2\leq C_2<C_1$. When $D_1> \a_1 C_3$, and the metering rate $D_2<\a_2 C_3$, then $q_1=\a_1 C_3$ and $q_2=D_2$. In this case, $q_3=\a_1 C_3+D_2$, and the utilization rate of link 3's capacity is 
\bqs
\frac{q_3}{\min\{ D_1+D_2,C_3\}}&=&\frac{\a_1 C_3+D_2}{\min\{ D_1+D_2,C_3\}}=\frac{\a_1 C_3+D_2}{C_3}\geq \a_1.
\eqs
That is, as much as $\a_2$ of link 3's capacity is wasted. 
In contrast, if we increase the metering rate of link 2, such that $D_2\geq \a_2 C_3$. Then \refe{invariant:constant} predicts $q_3=C_3$.
That is, low metering rates could cause a lower utilization rate of link 3's capacity. This property of the merge model \refe{invariant:constant} is counter-intuitive, and the merge model is not physical.

\section{Numerical examples}
In this section, we numerically solve various merge models and demonstrate the validity of our analytical results. 
Here, both links 1 and 3 are two-lane mainline freeways with a corresponding normalized maximum sensitivity fundamental diagram \citep{delcastillo1995fd_empirical} is ($\r\in[0,2]$)
\bqs
Q(\r)&=& \r \left\{ 1-\exp\left[1-\exp\left( \frac14(\frac{2}{\r}-1)\right)\right]\right\}.
\eqs
Link 2 is a one-lane on-ramp with a fundamental diagram as ($\r\in[0,1]$)
\bqs
Q(\r)&=& \frac 12 \r \left\{ 1-\exp\left[1-\exp\left( \frac14(\frac{1}{\r}-1)\right)\right]\right\}.
\eqs
Note that here the free flow speed on the on-ramp is half of that on the mainline freeway, which is 1. Thus we have the capacities $C_1=C_3=4C_2=0.3365$ and the corresponding critical densities $\r_{c1}=\r_{c3}=2\r_{c2}=0.4876$. The length of all three links is the same as $L=10$, and the simulation time is $T=360$.

In the following numerical examples, we discretize each link into $M$ cells and divide the simulation time interval $T$ into $N$ steps. The time step $\dt=T/N$ and the cell size $\dx=L/M$, with $\dt=0.9 \dx$, satisfy the CFL condition \citep{courant1928CFL}
\bqs
v_f \frac{\dt}{\dx}=\frac{\dt}{\dx} =0.9 \leq 1.
\eqs
Then we use the following finite difference equation  for link $i=1,2,3$:
\bqs
\r_{i,m}^{n+1}&=&\r_{i,m}^n+\frac{\dt}{\dx}(q_{i,m-1/2}^n-q_{i,m+1/2}^n),
\eqs 
where $\r_{i,m}^n$ is the average density in cell $m$ of link $i$ at time step $n$, and the boundary fluxes $q_{i,m-1/2}^n$ are determined by supply-demand methods. For example, for upstream links $i=$1 and 2, the in-fluxes are
\bqs
q_{i,m-1/2}^n&=&\min\{D_{i,m-1}^n,S_{i,m}^n\}, \quad m=1,\cdots, M,
\eqs
where $D_{i,m-1}^n$ is the demand of cell $m-1$ on link $i$, $S_{i,m}^n$ the supply of cell $m$, and $D_{i,0}^n$ the demand of link $i$. 
For link 3, the out-fluxes are
\bqs
q_{3,m+1/2}^n&=&\min\{D_{3,m}^n,S_{3,m+1}^n\}, \quad m=1,\cdots, M,
\eqs
where $S_{3,M+1}^n$ is the supply of the destination.
Then the out-fluxes of the upstream links and the in-flux of the downstream link are determined by merge models, which are discrete versions of \refe{sdentropy}:
\bqs
q_{i,M+1/2}^n&=&F_i(D_{1,M}^n, D_{2,M}^n, S_{3,1}^n),\\
q_{3,1/2}^n&=&q_{1,M+1/2}^n+q_{2,M+1/2}^n.
\eqs

In our numerical studies, we only consider the fair merge model \refe{sdentropy:fair} and its invariant counterpart \refe{sdentropy:priority}. That is, in the fair merge model, we have
\bqs
q_{i,M+1/2}^n&=&\min\{D_{i,M}^n,\frac{D_{i,M}^n}{D_{1,M}^n+ D_{2,M}^n} S_{3,1}^n\},\\
q_{3,1/2}^n&=&\min\{D_{1,M}^n+ D_{2,M}^n, S_{3,1}^n\}.
\eqs
In the invariant fair merge model, we have ($i,j=1, 2$ and $i\neq j$)
\bqs
q_{i,M+1/2}^n&=&\min\{D_{i,M}^n,\max\{S_{3,1}^n- D_{j,M}^n, \a_i S_{3,1}^n \} \},\\
q_{3,1/2}^n&=&\min\{D_{1,M}^n+ D_{2,M}^n, S_{3,1}^n\},
\eqs
where $\a_i=C_i/(C_1+C_2)$.

\subsection{Kinematic waves, stationary states, and interior states in the fair merge model}
In this subsection, we study numerical solutions of the fair merge model in \refe{sdentropy:fair}. Initially, links 1 and 3 carry UC flows with $\r_1=\r_3=0.35$, and at $t=0$ a traffic stream on link 2 starts to merge into link 3 with $\r_2=0.1$. That is, the initial conditions in supply-demand space is 
$U_1=U_3=(0.3131, 0.3365)$ and $U_2=(0.0500, 0.0841)$. 
Here we use the Neumann boundary condition in supply and demand \citep{collela2004_pde}: $D_{1,0}^n=D_{1,1}^n$, $D_{2,0}^n=D_{2,1}^n$, and $S_{3,M+1}^n=S_{3,M}^n$. Therefore, we have a Riemann problem here.

In this case, $D_1+D_2=0.3631 > S_{3}=0.3365$ and $D_2 \leq \frac{C_2}{C_1+C_2}S_3=0.0673$. Thus according to Theorem \ref{thm:fairmerg2}, we have the following stationary and interior states $U_{3}^+=U_3(0^+,t)=(C_{3}, S_{3})=(C_1,C_1)$, $U^-_2=(D_2, C_2)$, $U_2(0^-,t)=(\frac{C_1}{S_3-D_2} D_2,C_2)$, and $U^-_1=U_1(0^-,t)=(C_1, S_3-D_2)$. 
 Furthermore, from the LWR model, there will be a back-traveling shock wave on link 1 connecting $U_1$ to $U_1^-$, no wave on link 2, and a rarefaction wave on link 3 connecting $U_3^+$ to $U_3$.

\bfg\bc
\includegraphics[width=4in]{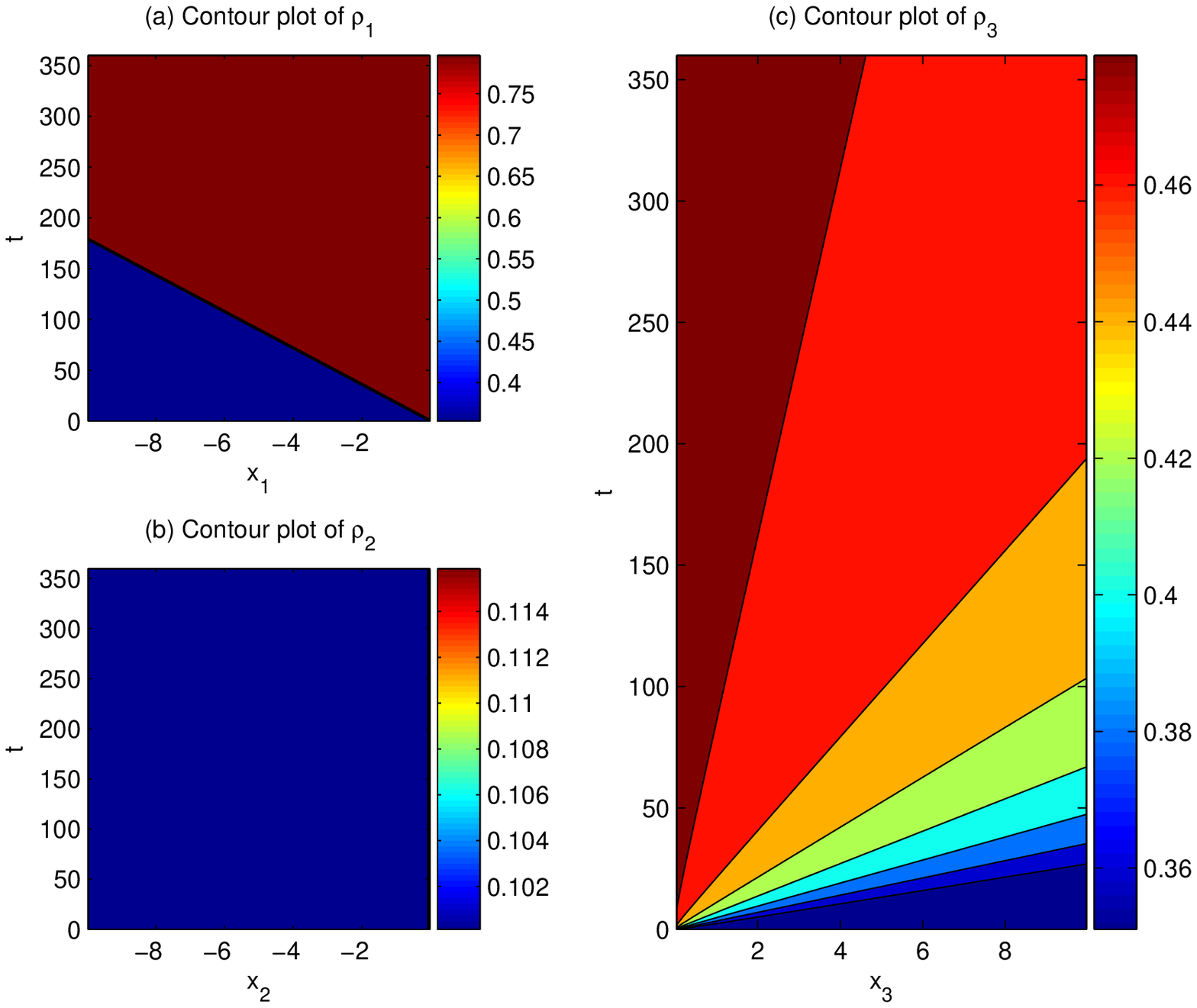} \caption{Solutions of the fair merge model \refe{sdentropy:fair}: $M=160$, $N=6400$.} \label{fairmerge20090219contour}
\ec\efg

In \reff{fairmerge20090219contour}, the solutions of $\r_1$, $\r_2$, and $\r_3$ are demonstrated with $M=160$ and $N=6400$. From the figures, we can clearly observe the predicted kinematic waves. 
In \reff{fairmerge20090219contour}(b) there is a very thin layer of higher densities in the last cell of the upstream link near the merging junction, which is caused by the interior state as predicted.
In addition, we can observe at $t=T$ the approximate asymptotic values: $U_1^-=U_1(0^-,t)=(0.3365, 0.2865)$, and $\r_1^-=\r_1(0^-,t)=0.8277$; $U_2^-=(0.05, 0.0841)$, $U_2(0^-,t)=(0.0587, 0.0841)$, $\r_2^-=0.1$, $\r_2(0^-,t)=0.1179$;  $U_3^+=(0.3365, 0.3365)$, and $\r_3^+=\r_{c3}=0.4874$. These numbers are all very close to the theoretical values and get closer if we reduce $\dx$ or increase $T$. That is, the results are consistent with theoretical results asymptotically.

\bfg\bc
\includegraphics[width=4in]{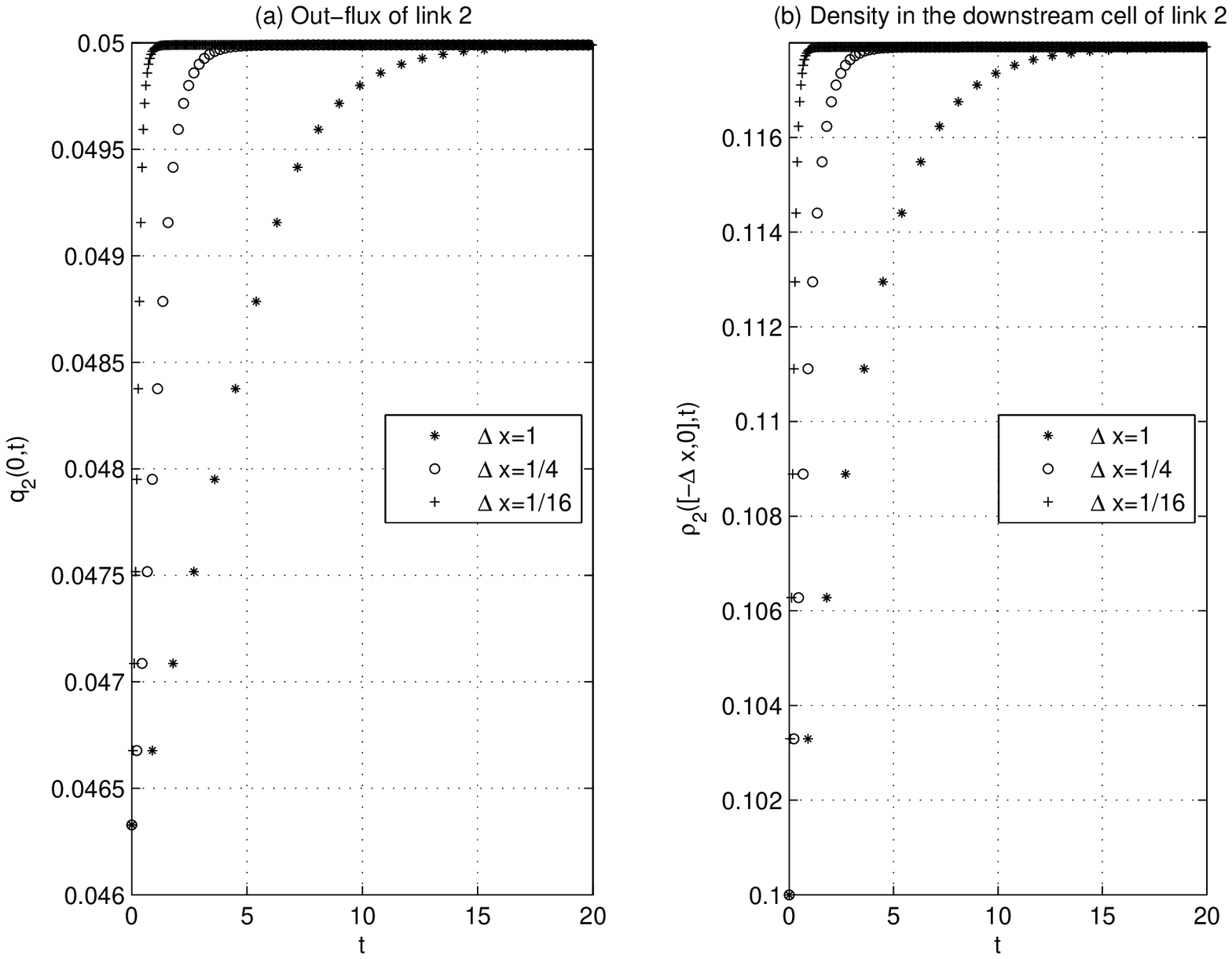} \caption{Evolution of the out-flux and the density in the downstream cell of link 2 for the fair merge model \refe{sdentropy:fair}} \label{fairmerge20090219interior}
\ec\efg

In \reff{fairmerge20090219interior}, we demonstrate the evolution of the out-flux and the density in the downstream cell of link 2 for three different cell sizes. From \reff{fairmerge20090219interior}(a) we can see that, initially, the out-flux of link 2 is
\bqs
\frac{D_2}{D_1+D_2}S_3=\frac{0.05}{0.05+0.3131} 0.3365 =0.0463,
\eqs
which is not the same but approaches the asymptotic out-flux $D_2=0.05$. Correspondingly the density in the downstream cell of link 2 approaches the interior state, as shown in \reff{fairmerge20090219interior}(b). Furthermore, as we decrease the cell size, the numerical results are closer to the theoretical ones at the same time instant. This figure shows that the fair merge model is not invariant, but approaches its invariant counterpart asymptotically. Note that the densities in any other cells of link 2 remain constant at 0.1.

\subsection{Comparison of non-invariant and invariant merge models}
In this subsection, we compare the numerical solutions of the fair merge model \refe{sdentropy:fair} with its invariant counterpart  \refe{sdentropy:priority2}. Initially, links 1 and 3 carry UC flows with $\r_1=\r_3=0.35$, and at $t=0$ a traffic stream on link 2 starts to merge into link 3 with $\r_2=0.1$. Different from the example in the preceding subsection, here we use the following boundary conditions: $D_{1,0}^n=D_{1,1}^n$, $D_{2,0}^n=0.05+0.03\sin(n\pi \dt /60)$, and $S_{3,M+1}^n=S_{3,M}^n$. Thus we have a periodic demand on link 2. 

\bfg\bc
\includegraphics[width=4in]{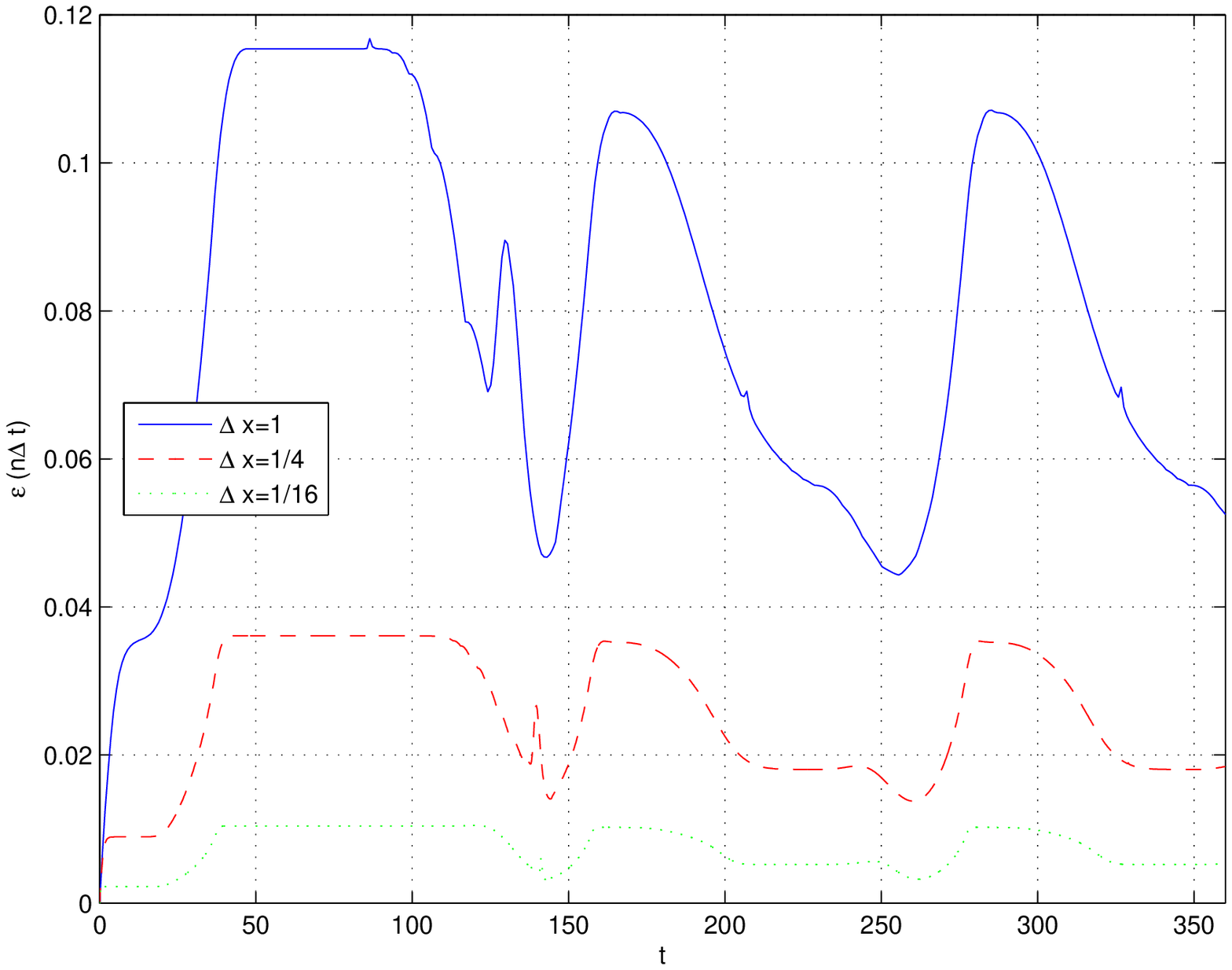}\caption{Difference in the solutions between the fair merge model \refe{sdentropy:fair} and its invariant counterpart  \refe{sdentropy:priority}} \label{merge20090219comparison}
\ec\efg

We use $\r_{i,m}^n$ for the discrete density from the fair merge model \refe{sdentropy:fair} and $\bar \r_{i,m}^n$ from its invariant counterpart  \refe{sdentropy:priority2}. Then we denote the difference between the two solutions by
\bqn
\epsilon (n \Delta t)&=& \sum_{i=1}^ 3 \sum_{m=1}^M |\r_{i,m}^n- \bar \r_{i,m}^n| \dx.
\eqn
In \reff{merge20090219comparison}, we can see that the difference decreases if we decreases the cell size. This clearly demonstrates that the fair merge model \refe{sdentropy:fair} converges to its invariant counterpart  \refe{sdentropy:priority}.

\section{{Conclusion}}
In this paper, we studied continuous kinematic wave models of merging traffic flow which are consistent with discrete CTM with various distribution schemes. In particular, we introduced the supply-demand diagram of traffic flow and proposed a solution framework for the Riemann problem of merging traffic flow. In the Riemann solutions, each link can have two new states, an interior state and a stationary state, and the kinematic waves on a link are determined by the initial state and the stationary state. We then derived admissible conditions for interior and stationary states and introduced various distribution schemes as entropy conditions defined in the interior states. Then we proved that stationary states and boundary fluxes exist and are unique for the Riemann problem for both the fair and constant distribution schemes. We also discussed two invariant merge models, in which the local and discrete flux is the same as the global and continuous flux. With numerical examples, we demonstrated the validity of the proposed analytical framework and that the fair merge model converges to its invariant counterpart when we decrease the cell size. Compared with existing discrete kinematic wave models (i.e. CTM with various distribution schemes) of merging traffic flow, the continuous models can provide analytical insights on kinematic waves arising at merging junctions; and, compared with existing analytical models, they are physically meaningful and consistent with existing CTM.

In this study, we introduced a new definition of invariant merge models, in which fluxes computed by discrete models should be the same as those by their continuous counterparts.
An important observation is that both the fair and constant merge models are not invariant. 
For example, for the fair merge model, at $t=0$ the local fluxes from \refe{sdentropy:fair} are
\bqs
q_i&=&\min\{1,\frac{S_3}{D_1+D_2}\} D_i,\quad i=1,2,
\eqs
which are different from \refe{sdentropy:fair2} when only one upstream is congested; i.e., when $D_1+D_2>S_3$ and $D_i \leq \frac{C_i}{C_1+C_2}S_3$. However, the results here suggest that the discrete fluxes converge to the continuous ones after a sufficient amount of time or at a given time but with decreasing time intervals. That is, the non-invariant discrete merge models do not provide incorrect solutions but just approximate solutions to the corresponding continuous merge models. Compared with invariant discrete merge models, the fair merging model has some unique merits; e.g., it can be easily extended to general junctions with multiple upstream links as shown in \citep{jin2003merge,jin2004network}. 
Note that, as demonstrated in \citep{jin2009sd}, invariant models can also yield interior states.	

This paper presents a systematic framework for solving kinematic waves arising from merging traffic in supply-demand space. We expect that boundary fluxes, stationary states, and kinematic waves for other distribution schemes can also be solved in this framework. 
For example, we can obtain the stationary states and kinematic waves for the following artificial merge model
\bqs
q_i&=&\min\{1,\frac{0.9 S_3(0^+,t)}{D_1(0^-,t)+D_2(0^-,t)}\} D_i(0^-,t),\quad i=1,2,
\eqs
in which only 90\% of the downstream supply can be utilized by the upstream traffic. 

The Riemann problem for a merge with three or more upstream links can be discussed in the same framework, but, due to the space limitations, this will be discussed in future studies. Generally, there are systematic lane-changing activities at merging junctions, and existing merge models based on supply-demand method cannot capture the impacts of lane-changes \citep{laval2005impacts}. In the future, it would be interesting to analyze the formation and dissipation of traffic queues at merging junctions when lane-changes are considered.

\section*{Acknowledgements}
I would like to thank Dr. Jorge A. Laval of Georgia Tech for his helpful comments on an earlier version of the paper. Constructive comments of two anonymous reviewers have been very helpful for improving the presentation of the paper. The views and results contained herein are the author's alone.

\pdfbookmark[1]{Appendix A}{appendixa}
\section*{Appendix A: Proof of Theorem \ref{thm:fairmerg2}}
{\em Proof}. From traffic conservation equations in \refe{trafficconservation} and admissible conditions of stationary states, we can see that 
\bqs
q_3\leq \min\{D_1+D_2, S_3 \}.
\eqs
We first demonstrate that it is not possible that $q_3 < \min\{D_1+D_2, S_3 \} \leq \min\{C_1+C_2, C_3\}$. Otherwise, from \refe{downstreamss} and \refe{downstreamis} we have $U_3(0^+,t)=U_3^+=(q_3,C_3)$ with $q_3<S_3$; Since $q(U_1^-)+q(U_2^-)=q_3<D_1+D_2$, then we have $q(U_i^-)<D_i$ for at least one upstream link, e.g., $q_1<D_1$. From \refe{upstreamss} and \refe{upstreamis} we have $U_1(0^-,t)=U_1^-=(C_1, q_1)$.
Then from the entropy condition in \refe{sdentropy:fair} we have
\bqs
q_3&=&\min\{C_1+D_2(0^-,t), C_3\},\\
q_1&=&\min\{1,\frac{C_3}{C_1+D_2(0^-,t)}\} C_1.
\eqs
Since $q_3<C_3$, from the first equation we have $q_3=C_1+D_2(0^-,t)<C_3$, and from the second equation we have $q_1=C_1$, which contradicts $q_1<D_1$. Therefore, 
\bqs
q_1+q_2=q_3 = \min\{D_1+D_2, S_3 \}.
\eqs
That is, the fair distribution scheme yields the optimal fluxes for any initial conditions.

\bi
\item [(1)] When $D_1+D_2<S_3$, we have $q_3=D_1+D_2<S_3$. From \refe{downstreamss} and \refe{downstreamis} we have $U_3(0^+,t)=U_3^+=(D_1+D_2,C_3)$. 
Since $q_1+q_2=D_1+D_2$ and $q_i\leq D_i$, we have $q_i= D_i$. From \refe{upstreamss} we have $U_i^-=(D_i,C_i)$. From \refe{upstreamis} we have $U_i(0^-,t)=(D_i(0^-,t),S_i(0^-,t))$ with $S_i(0^-,t)\geq D_i^-=D_i$. From \refe{sdentropy:fair} we have
\bqs
q_3&=&\min\{D_1(0^-,t)+D_2(0^-,t), C_3\}=D_1+D_2<S_3\leq C_3,\\
q_i&=&\min\{1,\frac{C_3}{D_1(0^-,t)+D_2(0^-,t)}\} D_i(0^-,t) =D_i.
\eqs
Thus, $D_i(0^-,t)=D_i\leq S_i(0^-,t)$. Then $U_i(0^-,t)=U_i^-=(D_i,C_i)$. In this case, there are no interior states on all links.

\item [(2)] When $D_1+D_2=S_3$, we have $q_3=S_3$. From \refe{downstreamss} and \refe{downstreamis}, we have $U_3^+=(C_3,S_3)$ and $U_3(0^+,t)=(D_3(0^+,t),S_3(0^+,t))$ with $D_3(0^+,t)\geq S_3^+=S_3$. 
Since $q_1+q_2=D_1+D_2$ and $q_i\leq D_i$, we have $q_i= D_i$. From \refe{upstreamss} and \refe{upstreamis}, we have $U_i^-=(D_i,C_i)$ and $U_i(0^-,t)=(D_i(0^-,t),S_i(0^-,t))$ with $S_i(0^-,t)\geq D_i^-=D_i$. 
From \refe{sdentropy:fair} we have
\bqs
q_3&=&\min\{D_1(0^-,t)+D_2(0^-,t), S_3(0^+,t)\}=D_1+D_2=S_3,\\
q_i&=&\min\{1,\frac{S_3(0^+,t)}{D_1(0^-,t)+D_2(0^-,t)}\} D_i(0^-,t) =D_i.
\eqs
We can have the following two scenarios. 
\bi
\item [(i)] If $D_1(0^-,t)+D_2(0^-,t)\geq S_3(0^+,t)=D_1+D_2=S_3 \leq D_3(0^+,t)$, then $U_3(0^+,t)=U_3^+=(C_3,S_3)$ and there is no interior state on link 3. Moreover, we have 
\bqs \frac{D_1+D_2}{D_1(0^-,t)+D_2(0^-,t)} D_i(0^-,t) =D_i,
\eqs
 which leads to $D_i(0^-,t)\leq D_i\leq S_i(0^-,t)$. From the assumption that $D_1(0^-,t)+D_2(0^-,t)\geq D_1+D_2$, we have $D_i(0^-,t) = D_i$. Further we have $U_i(0^-,t)=U_i^-=(D_i,C_i)$, and there are no interior states on links 1 or 2.
\item [(ii)] If $S_3(0^+,t)>D_1(0^-,t)+D_2(0^-,t)=D_1+D_2=S_3$, $D_i(0^-,t) =D_i$. Thus $U_i(0^-,t)=U_i^-=(D_i,C_i)$, and there are no interior states on links 1 or 2. Moreover, $U_3(0^+,t)$ satisfies $S_3(0^+,t)>S_3$ and $D_3(0^+,t)\geq S_3$. Thus there can be multiple interior states on link 3 when $S_3<C_3$.
\ei

\item [(3,4)]
When $D_1+D_2>S_3$, for upstream links, at least one of the stationary states is SOC. Otherwise, from \refe{upstreamss} we have $U_i^-=(D_i,C_i)$, and $q_1+q_2=D_1+D_2>S_3$, which is impossible. In addition, we have $q_3=S_3<D_1+D_2$. From \refe{downstreamss} we have $U_3^+=(C_3,S_3)$. From \refe{downstreamis} we have $U_3(0^+,t)=(D_3(0^+,t),S_3(0^+,t))$ with $D_3(0^+,t)\geq S_3^+=S_3$. 	 
From \refe{sdentropy:fair} we have
\bqs
q_3&=&\min\{D_1(0^-,t)+D_2(0^-,t), S_3(0^+,t)\}=S_3<D_1+D_2,\\
q_i&=&\min\{1,\frac{S_3(0^+,t)}{D_1(0^-,t)+D_2(0^-,t)}\} D_i(0^-,t).
\eqs
If $D_1(0^-,t)+D_2(0^-,t)\leq S_3(0^+,t)$, then $D_1(0^-,t)+D_2(0^-,t)=S_3<D_1+D_2$ and $q_i=D_i(0^-,t)$. This is not possible for the SOC stationary state $U_i^-=U_i(0^-,t)=(C_i,q_i)$ with $q_i<D_i\leq C_i$. Thus $S_3(0^+,t)<D_1(0^-,t)+D_2(0^-,t)$, $S_3(0^+,t)=S_3<D_1+D_2$, and $U_3(0^+,t)=U_3^+=(C_3,S_3)$. Hence for both upstream links
\bqs
q_i&=&\frac{S_3}{D_1(0^-,t)+D_2(0^-,t)}D_i(0^-,t).
\eqs

\bi
\item [(3)] When $D_i> \frac{C_i}{C_1+C_2}S_3$ ($i=1,2$), stationary states on both links 1 and 2 are SOC with $U_i^-=U_i(0^-,t)=(C_i,q_i)$ with $q_i<D_i$. Otherwise, we assume that link 1 is SOC with $U_1(0^-,t)=U_1^-=(C_1,q_1)$ and link 2 is UC with $U_2^-=(D_2,C_2)$. Then 
\bqs
D_2&=&\frac{S_3}{C_1+D_2(0^-,t)}D_2(0^-,t) \leq \frac{S_3}{C_1+C_2}C_2<D_2,
\eqs
which is impossible.
From \refe{sdentropy:fair}, we have 
\bqs
q_i=\frac{S_3}{C_1+C_2} C_i,
\eqs
and $U_i(0^-,t)=U_i^-=(C_i,q_i)$.

\item[(4)]
When $D_1+D_2 > S_{3}$ and $D_i \leq \frac{C_i}{C_1+C_2}S_3$ ($i,j=1$ or 2 and $i\neq j$), we can show that stationary states on links $j$ and $i$ are SOC and UC respectively with $U_j^-=U_j(0^-,t)=(C_j,q_j)$ with $q_j<D_j$, $U_i^-=(D_i,C_i)$, and $S_i(0^-,t)\geq D_i$. Otherwise, $U_i(0^-,t)=U_i^-=(C_i,q_i)$ with $q_i<D_i$, and
\bqs
q_i&=&\frac{S_3}{C_i+D_j(0^-,t)} C_i \geq \frac{C_i}{C_1+C_2}S_3 \geq D_i,
\eqs
which is impossible. Since at least one of the upstream links has SOC stationary state, the stationary states on links $i$ and $j$ are UC and SOC respectively.
 From \refe{sdentropy:fair}, we have a unique interior state on link $i$, $U_i(0^-,t)=(\frac{D_i}{S_3-D_i}C_j, C_i)$, and $q_j=S_3-D_i$.
\ei
\ei
For the four cases, it is straightforward to show that \refe{sdentropy:fair2} always holds.
\eop

\pdfbookmark[1]{Appendix B}{appendixb}
\section*{Appendix B: Proof of Theorem \ref{thm:constmerg2}}
{\em Proof}. 
\bi
\item [(1)]
When $D_1+D_2 < S_{3}$ and $D_i\leq \a_i C_3$ ($i=1,2$), then $q_3<S_3$, and $U_3^+=U_3(0^+,t)=(q_3,C_3)$. Assuming that $q_i<D_i\leq C_i$, then $U_i^-=U_i(0^-,t)=(C_i,q_i)$. From \refe{sdentropy:constant}, we have $q_i=\min\{C_i, \a_i C_3\}=\a_i C_3<D_i$, which contradicts $D_i\leq \a_i C_3$. Therefore, $q_i=D_i$ and $q_3=D_1+D_2$.

\item [(2)]
When $D_i >\a_i C_3$ and $D_j<S_3-\a_i C_3$ ($i,j=$ 1 or 2 and $i\neq j$), we first prove that $q_i<D_i$ and then $q_j=D_j$.
\bi
\item[(i)] If $q_i=D_i$, then from \refe{sdentropy:constant}, we have $\a_i C_3 < D_i=\min\{D_i(0^-,t), \a_i S_3(0^+,t)\} \leq \a_i S_3(0^+,t)$, which leads to $C_3< S_3(0^+,t)$. This contradicts $C_3 \geq S_3(0^+,t)$. Thus, $q_i<D_i$, and $U_i^-=U_i(0^-,t)=(C_i,q_i)$. 

\item[(ii)] If $q_j<D_j\leq C_j$, then $U_j^-=U_j(0^-,t)=(C_j,q_j)$. From \refe{sdentropy:constant}, we have $q_j=\min\{C_j, \a_j S_3(0^+,t)\}=\a_j S_3(0^+,t)<D_j<S_3-\a_i C_3\leq \a_j C_3$. Thus $S_3(0^+,t)<C_3$.
Since $q_i = \a_i S_3(0^+,t)$, we have $q_3=q_i+q_j<\a_i S_3(0^+,t)+S_3-\a_i C_3<S_3$, which leads to $U_3^+=U_3(0^+,t)=(q_3,C_3)$. This  contradicts $S_3(0^+,t)<C_3$. Therefore, $q_j=D_j$. 
\ei
From \refe{sdentropy:constant}, we have $q_i=\a_i S_3(0^+,t)$. Thus $q_3=q_i+q_j <S_3-\a_i C_3+\a_i S_3(0^+,t)\leq S_3$. Then, $U_3^+=U_3(0^+,t)=(q_3,C_3)$, $q_i=\a_i C_3$, $q_3=\a_i C_3+D_j$.

\item[(3)]
When $D_1+D_2\geq S_3$, $S_3-\a_j C_3 \leq D_i \leq \a_i S_3$ ($i,j=$1 or 2 and $i\neq j$), we first prove that $q_3=S_3$ and then that $q_i=D_i$. Therefore $q_i=D_i$, and $q_j=S_3-D_i$.

\bi
\item [(i)]
If $q_3<S_3$, then $U_3^+=U_3(0^+,t)=(q_3,C_3)$. We first prove that at least one upstream stationary state is SOC and then that none of the upstream stationary states can be SOC. Therefore, $q_3=S_3$.
\bi
\item [(a)]
 If none of the upstream stationary states are SOC, then $q_i=D_i$ ($i=1,2$), and $q_3=q_1+q_2=D_1+D_2\geq S_3$, which contradicts $q_3<S_3$.
\item [(b)]
Assuming that $q_i<D_i$, then $U_i^-=U_i(0^-,t)=(C_i,q_i)$. From \refe{sdentropy:constant}, we have $q_i=\min\{C_i, \a_i C_3\}=\a_i C_3<D_i$. This is not possible, since $D_i\leq \a_i S_3$. Thus, $q_i=D_i$.
\item [(c)]
Assuming that $q_j<D_j$, then $U_j^-=U_j(0^-,t)=(C_j,q_j)$. From \refe{sdentropy:constant}, we have $q_j=\min\{C_j, \a_j C_3\}=\a_j C_3$. Then $D_i=q_i=q_3-\a_j C_3 <S_3 -\a_j C_3$, which contradicts $D_i\geq S_3-\a_jC_3$.
\ei

\item [(ii)]
Now assume $q_i<D_i$, then $U_i^-=U_i(0^-,t)=(C_i,q_i)$. From \refe{sdentropy:constant}, we have $q_i=\min\{C_i, \a_i S_3(0^+,t)\}=\a_i S_3(0^+,t)<D_i\leq \a_i S_3$. Therefore, $q_j\leq \a_j S_3(0^+,t) < \a_j S_3$, which leads to $q_i+q_j<S_3$. This is not possible, since $q_i+q_j=S_3$.
\ei

\item[(4)]
When $D_i \geq \alpha_i S_3$, we first prove that $q_3=S_3$ and then that $q_i \geq \a_i S_3$ ($i=1,2$). Therefore, $q_i = \a_i S_3$ ($i=1,2$), and $q_3=S_3$.
\bi
\item [(i)] If $q_3<S_3$, then $U_3^+=U_3(0^+,t)=(q_3,C_3)$. We first prove that at least one upstream stationary state is SOC and then that none of the upstream stationary states can be SOC. Therefore $q_1+q_2=q_3=S_3$.

\bi
\item [(a)] If none of the upstream stationary states are SOC, then $q_i=D_i$ ($i=1,2$), and $q_3=D_1+D_2\geq S_3$, which contradicts $q_3<S_3$.
\item [(b)] Assuming that $q_i<D_i$, then  $U_i^-=U_i(0^-,t)=(C_i,q_i)$. From \refe{sdentropy:constant}, we have $q_i=\min\{C_i, \a_i C_3\}=\a_i C_3$. Then $q_j=q_3-q_i<S_3-\a_i C_3 \leq \a_j S_3 \leq D_j$. Thus, $U_j^-=U_j(0^-,t)=(C_j,q_j)$, and from \refe{sdentropy:constant} we have $q_j=\min\{C_j, \a_j C_3\}=\a_j C_3>\a_j S_3$, which contradicts $q_j <\a_j S_3$. Thus, $U_i^-$ is UC and $q_i=D_i$.
\item [(c)] Similarly we can prove that $U_j^-$ is UC and $q_j=D_j$.
\ei
 
\item [(ii)]
If $q_i<\a_i S_3\leq D_i\leq C_i$, then $U_i^-=U_i(0^-,t)=(C_i,q_i)$. From \refe{sdentropy:constant}, we have $q_i=\min\{C_i,\a_i S_3(0^+,t)\}=\a_i S_3(0^+,t) < \a_i S_3$. Therefore, $q_j\leq \a_j S_3(0^+,t)<\a_jS_3$, which leads to $q_i+q_j<S_3$. This is not possible, since $q_i+q_j=S_3$.
\ei
\ei
For the four cases, it is straightforward to show that \refe{sdentropy:constant2} always holds.
\eop

\pdfbookmark[1]{Appendix C}{appendixc}
\section*{Appendix C: Proof of Theorem \ref{thm:invariant1}}
{\em Proof}. 

First \refe{sdentropy:priority} implies that
\bqs
q_3&=&\min\{D_1(0^-,t)+D_2(0^-,t), S_3(0^+,t) \},
\eqs
which can be shown for three cases: (i) $D_1(0^-,t)+D_2(0^-,t)<S_3(0^+,t)$, (ii) $D_i(0^-,t)\geq \a_iS_3(0^+,t)$, and (iii) $D_1(0^-,t)+D_2(0^-,t)\geq S_3(0^+,t)$ and $D_i(0^-,t)\leq \a_iS_3(0^+,t)$.
\bi
\item [(1)] When $D_1+D_2<S_3$, $q_3=q_1+q_2\leq D_1+D_2<S_3\leq C_3$. Thus the downstream stationary state is SUC with $U_3^+=U_3(0^+,t)=(q_3,C_3)$. In the following, we prove that $q_i=D_i$, which is consistent with \refe{sdentropy:priority2}. Therefore \refe{sdentropy:priority2} is correct in this case.
\bi
\item[(i)] Assuming that $q_i<D_i\leq C_i$, then the stationary state on link $i$ is SOC with $U_i^-=U_i(0^-,t)=(C_i,q_i)$. From \refe{sdentropy:priority}, we have
\bqs
q_i&=&\min\{C_i, \max\{C_3-D_j(0^-,t), \a_i C_3\}\}=\max\{C_3-D_j(0^-,t), \a_i C_3\}<C_i,\\
q_j&=&\min\{D_j(0^-,t), \max\{C_3-C_i, \alpha_j C_3\}\}.
\eqs
We show that the two equations have no solutions for either $\a_j C_3 \leq D_j(0^-,t)$ or $\a_j C_3 > D_j(0^-,t)$. Thus $q_i=D_i$.

\bi
\item [(a)] When $\a_j C_3 \leq D_j(0^-,t)$, we have $\a_i C_3 \geq C_3-D_j(0^-,t)$. From the first equation we have $q_i=\a_i C_3$. From the second equation we have $q_j= D_j(0^-,t)\geq \a_j C_3$ or $q_j=\max\{C_3-C_i, \alpha_j C_3\}\geq \a_j C_3$. Thus $q_i+q_j\geq C_3\geq S_3$, which contradicts $q_3<S_3$.
\item [(b)] When $\a_j C_3 > D_j(0^-,t)$, we have $\a_i C_3 < C_3-D_j(0^-,t)$. From the first equation we have $q_i= C_3-D_j(0^-,t)$. From the second equation we have $q_j=D_j(0^-,t)$. Thus $q_i+q_j=C_3$, which contradicts $q_3<S_3$.
\ei
\ei

\item [(2)] When $D_i \geq \a_i S_3$, $S_3-D_j \leq \a_i S_3$. In the following we show that $q_3=S_3$ and $q_i=\a_i S_3$, which is consistent with \refe{sdentropy:priority2}.
\bi
\item [(i)] If $q_3<S_3$, then the stationary state on link 3 is SUC with $U_3^+=U_3(0^+,t)=(q_3,C_3)$. Also at least one of the upstream stationary states is SOC, since, otherwise, $q_1+q_2=D_1+D_2\geq S_3$. Here we assume $U_i^-=U_i(0^-,t)=(C_i,q_i)$.
From  \refe{sdentropy:priority} we have
\bqs
q_i&=&\min\{C_i, \max\{C_3-D_j(0^-,t), \alpha_i C_3\}\}=\max\{C_3-D_j(0^-,t), \alpha_i C_3\},\\
q_j&=&\min\{D_j(0^-,t), \max\{C_3-C_i, \alpha_j C_3\}\}.
\eqs 
We show that the two equations have no solutions for either $C_3-D_j(0^-,t) \geq \a_i C_3$ or $C_3-D_j(0^-,t) < \a_i C_3$. Thus $q_3=S_3$.
\bi
\item [(a)] If $C_3-D_j(0^-,t) \geq \a_i C_3$, $D_j(0^-,t) \leq \a_j C_3$. From the first equation we have $q_i= C_3-D_j(0^-,t)$. From the second equation we have $q_j=D_j(0^-,t)$. Thus $q_i+q_j=C_3$, which contradicts $q_3<S_3\leq C_3$.
\item [(b)] If $C_3-D_j(0^-,t) < \a_i C_3$, $D_j(0^-,t) > \a_j C_3$. From the first equation we have $q_i= \a_i C_3$. From the second equation we have $q_j=D_j(0^-,t) > \a_j C_3$ or $q_j=\max\{C_3-C_i, \alpha_j C_3\}\geq \a_j C_3$. Thus $q_i+q_j\geq C_3$, which contradicts $q_3<S_3\leq C_3$.
\ei
\item [(ii)] If $q_i<\a_i S_3\leq D_i\leq C_i$ for any $i=1,2$, then $U_i^-=U_i(0^-,t)=(C_i,q_i)$. From  \refe{sdentropy:priority} we have
\bqs
q_i&=&\max\{S_3(0^+,t)-D_j(0^-,t), \alpha_i S_3(0^+,t)\}<C_i,\\
q_j&=&\min\{D_j(0^-,t), \max\{S_3(0^+,t)-C_i, \alpha_j S_3(0^+,t)\}\}.
\eqs
The first equation implies that $\alpha_i S_3(0^+,t)< \a_i S_3$; i.e., $S_3(0^+,t)<S_3$. In addition, $S_3(0^+,t)-D_j(0^-,t)<\a_i S_3$. Thus, $S_3(0^+,t)-C_i<S_3-C_i<S_3-\a_i S_3=\a_j S_3$, and $\max\{S_3(0^+,t)-C_i, \alpha_j S_3(0^+,t)\}<\a_j S_3$. From the second equation we have $q_j<\a_j S_3$. Thus $q_i+q_j<S_3$, which contradicts $q_i+q_j=S_3$. Thus $q_i\geq \a_i S_3$ for $i=1,2$. Since $q_i+q_j=S_3$, $q_i=\a_i S_3$.
\ei

\item[(3)] When $D_i+D_j \geq S_3$ and $D_i\leq \a_i S_3$ for $i,j=1$ or 2 and $i\neq j$. In the following we show that $q_3=S_3$ and $q_i=D_i$, which is consistent with \refe{sdentropy:priority2}.
\bi
\item [(i)] If $q_3<S_3$, then the stationary state on link 3 is SUC with $U_3^+=U_3(0^+,t)=(q_3,C_3)$. We first prove that at least one upstream stationary state is SOC and then that none of the upstream stationary states can be SOC. Therefore, $q_3=S_3$.

\bi
\item [(a)] If none of the upstream stationary states are SOC, then $q_1+q_2=D_1+D_2\geq S_3$, which contradicts $q_3<S_3$. Thus, at least one of the upstream stationary states is SOC.
\item [(b)] Assuming that $q_i<D_i$, then  $U_i^-=U_i(0^-,t)=(C_i,q_i)$.
From  \refe{sdentropy:priority} we have 
\bqs
q_i&=&\min\{C_i, \max\{C_3-D_j(0^-,t), \alpha_i C_3\}\}=\max\{C_3-D_j(0^-,t), \alpha_i C_3\}<D_i,
\eqs 
which is not possible, since $D_i\leq \a_i S_3$. Thus $q_i=D_i$.
\item [(c)] Assuming that $q_j<D_j$, then  $U_j^-=U_j(0^-,t)=(C_j,q_j)$. Since $q_3=\min\{D_i(0^-,t)+C_j, C_3 \}<S_3$, we have $q_3=D_i(0^-,t)+C_j<S_3$. 
From  \refe{sdentropy:priority} we have $D_i=q_i\leq D_i(0^-,t)$. Thus $D_i+D_j\leq D_i(0^-,t)+C_j <S_3$, which contradicts $D_i+D_j\geq S_3$.
\ei

\item [(ii)] If $q_i<D_i$, then $U_i^-=U_i(0^-,t)=(C_i,q_i)$. From \refe{sdentropy:priority}, we have 
\bqs
q_i&=&\max\{S_3(0^+,t)-D_j(0^-,t), \alpha_i S_3(0^+,t)\}<D_i\leq \a_i S_3,\\
q_j&=&\min\{D_j(0^-,t), \max\{ S_3(0^+,t)-C_i, \a_j S_3(0^+,t)\} \}.
\eqs
From the first equation we have that $S_3(0^+,t)<S_3$. We show that the two equations have no solutions for either $S_3(0^+,t)-D_j(0^-,t) \geq \alpha_i S_3(0^+,t)$ or $S_3(0^+,t)-D_j(0^-,t) < \alpha_i S_3(0^+,t)$. Therefore $q_i=D_i$.
\bi
\item [(a)] When $S_3(0^+,t)-D_j(0^-,t) \geq \alpha_i S_3(0^+,t)$, we have $D_j(0^-,t) \leq \a_j S_3(0^+,t)$. Thus $q_i=S_3(0^+,t)-D_j(0^-,t)$ and  $q_j=D_j(0^-,t)$. Then $q_i+q_j=S_3(0^+,t)<S_3$, which contradicts $q_i+q_j=S_3$.
\item [(b)] When $S_3(0^+,t)-D_j(0^-,t) < \alpha_i S_3(0^+,t)$, we have $q_i=\a_i S_3(0^+,t)$ and $S_3(0^+,t)-C_i < S_3(0^+,t)-q_i=\a_j S_3(0^+,t)$. Thus $q_j\leq \a_j S_3(0^+,t)$, and $q_i+q_j\leq S_3(0^+,t)<S_3$, which contradicts $q_3=S_3$. 
\ei

\ei

\ei

\eop

\end {document}